\documentclass[a4paper,11pt, english,twoside,centertags,leqno]{article}

	\usepackage{amsmath, amssymb,latexsym,amsxtra,amscd}
	\usepackage{stmaryrd} 
	\usepackage[all]{xy}
	\usepackage{theorem}
	\usepackage{rotating}
	
	\usepackage{verbatim}

	\usepackage{epic}
	\usepackage{graphicx}
	\usepackage{calc} 
	\graphicspath{{Figures/}}
	\usepackage{subcaption}
	\usepackage[usenames,dvipsnames,svgnames]{xcolor}
	\usepackage[shortlabels]{enumitem} 
		
	\theoremstyle{plain}

	{\theorembodyfont{\slshape}
	
	        \newtheorem{thm}{Theorem}[section]
	        
	        \newtheorem{lem}[thm]{Lemma}
	        \newtheorem{prop}[thm]{Proposition}
	        
	}
	
	{\theorembodyfont{\rmfamily}

	        \newtheorem{defn}[thm]{Definition}
	        
	        \newtheorem{rem}[thm]{Remark}

	}
	
	\def\l@section{\@tocline{1}{0pt}{0pc}{}{}}
	\def\l@subsection{\@tocline{2}{0pt}{1pc}{}{}}
	\def\l@subsubsection{\@tocline{3}{0pt}{2pc}{}{}}
	
	\AtBeginDocument{%
	}

	\def\cO{\mathcal{O}}

	\def\AA{\mathbb{A}}
	\def\CC{\mathbb{C}}
	\def\NN{\mathbb{N}}
	\def\pr{\mathbb{P}}
	
	\def\RR{\mathbb{R}}
	\def\ZZ{\mathbb{Z}}

	\let\oldvee\vee
	\renewcommand{\vee}{{\scriptscriptstyle\oldvee}}

	\let\epsilon\varepsilon
	\def\loccit{loc.\kern3pt cit.{}\xspace}
	\def\cf{cf.\kern.3em}
	\def\eg{e.g.\kern.3em}
	
	\def\resp{\text{resp.}\kern.3em}
	\let\setminus\smallsetminus
	\let\leq\leqslant

	\let\tilde\widetilde
	\let\bar\overline

	\newcommand{\ptbl}{.\kern .2em }

	\DeclareMathOperator{\id}{Id}

	\DeclareMathOperator{\GL}{GL}
	
	\DeclareMathOperator{\Sp}{Sp}

\newcommand{\E}{\mathcal{E}}
	\renewcommand{\pr}{{\frak p}}

	\newcommand{\g}{{\bf g}}

	\newcommand{\cx}{\tilde{x}}

	\newcommand{\To}{\;\longrightarrow\;}
	
	\newcommand{\liso}{\;\stackrel{\sim}{\longrightarrow}\;}

	\renewcommand{\pr}{{{\rm pr}}}

	\newcommand{\Perv}{{\rm Perv}}

	\newcommand{\End}{{\rm  End}}

	\newcommand{\cL}{{\mathcal{L}}}

	\renewcommand{\Im}{{\rm Im}}

	\newcommand{\PP}{{\mathbb{P}}}
	\newcommand{\cK}{{\mathcal{K}}}
	\newcommand{\cP}{{\mathcal{P}}}
	\newcommand{\cB}{{\mathcal{B}}}

	\def\to{\mathchoice{\longrightarrow}{\rightarrow}{\rightarrow}{\rightarrow}}

	\def\To#1{\mathchoice{\xrightarrow{\textstyle\kern4pt#1\kern3pt}}{\stackrel{#1}{\longrightarrow}}{}{}}

	\def\isoTo#1{\xrightarrow[\sim]{\textstyle\kern4pt#1\kern3pt}}
	
	\let\oldbigcup\bigcup
	\let\oldbigcap\bigcap
	\let\oldbigoplus\bigoplus
	\let\oldcoprod\coprod
	\def\bigcup{\mathop{\textstyle\oldbigcup}\displaylimits}
	\def\bigcap{\mathop{\textstyle\oldbigcap}\displaylimits}
	\def\bigoplus{\mathop{\textstyle\oldbigoplus}\displaylimits}
	\def\coprod{\mathop{\textstyle\oldcoprod}\displaylimits}

	\newcommand{\RedefinitSymbole}[1]{%
	\expandafter\let\csname old\string#1\endcsname=#1
	\let#1=\relax
	\newcommand{#1}{\csname old\string#1\endcsname\,}%
	}
	\RedefinitSymbole{\forall} \RedefinitSymbole{\exists}
	
	\usepackage[english]{babel}
	\usepackage[textwidth=15cm,textheight=21cm,heightrounded,hcentering]{geometry}
	\usepackage{titlesec}
	\titleformat{\section}[block]{\large\bfseries}{\thesection.}{0.5em}{}	
	\titleformat{\subsection}[runin]{\normalfont\bfseries}{\thesubsection.}{0.5em}{}[.]
	\usepackage[utf8]{inputenc}
	\usepackage[T1]{fontenc}
	
	\numberwithin{equation}{section}

	\usepackage[T1]{fontenc}
	\newcommand{\proof}{\noindent{\bf Proof:\ }}
	\newcommand{\Endproof}{\hspace*{\fill} $\Box$ \vspace{1ex} \noindent }
	
	\newcommand{\cR}{\mathcal{R}}

\usepackage{listings}
	\usepackage{hyperref} 

	\title{Monodromy of the Radon transform}

	\author{Meirav \textsc{Amram}\footnote{M.~Amram: Mathematics Department, Shamoon College of Engineering, 77245 Ashdod, Israel;	{meiravt@sce.ac.il}}\\
Michael \textsc{Dettweiler}\footnote{M.~Dettweiler: Department of Mathematics,
	University of Bayreuth, 95440 Bayreuth, Germany;
	{michael.dettweiler@uni-bayreuth.de}}\\
Stefan \textsc{Reiter}\footnote{S.~Reiter: Department of Mathematics,
	University of Bayreuth, 95440 Bayreuth, Germany;
	{stefan.reiter@uni-bayreuth.de}}}

\begin{document}
	\newgeometry{height=26cm,includefoot, includehead} 
	\maketitle

	\begin{abstract} We algorithmically
	determine the monodromy of the local system on the smooth part of  the Radon transformation of a generic simple perverse sheaf on the projective plane.
		\end{abstract}
		
	\tableofcontents
	
%


	\section*{Introduction}
	\addcontentsline{toc}{section}{Introduction}
	The classical Radon transform is an integral transformation which sends
	 a suitably integrable
	function $f$ on the projective plane $\PP^2(\RR)$
	to the  function
	$$\cR(f)\,:\,  \tilde{\PP}^2(\RR)\to \RR, \,\,\, \tilde{x}\mapsto \int_{L_{\tilde{x}}}f\,, $$ on the dual plane $\tilde{\PP}^2(\RR),$ where $L_{\tilde{x}}\subseteq \PP^2(\RR)$ is the line corresponding to an element $\tilde{x}$ in the dual plane \cite{Radon}. \\
	
		The algebraic-geometric version of  the classical Radon transformation
	was  introduced by Brylinsky  \cite{Brylinsky}.
	Following Grothendieck's fonctions-faisceaux philosophy (cf.~\cite{Brylinsky}, \cite{Laumon}),
	the basic idea is to replace functions by objects in the derived category
	$D^b_c(\PP^2)$ and integration of functions by direct image functors between derived categories.
		This leads to a functor (the {\it Radon transform})
	$$\cR: D^b_c(\PP^2) \to D^b_c(\tilde{\PP}^2),\,\, K\mapsto Rp_{2*}p_1^*(K)[1],$$ where
	$p_1,p_2$ are the two canonical projections of the incidence subvariety $H\subseteq
	\PP^2\times \tilde{\PP}^2$
	(see \cite{Brylinsky} and Sec.~\ref{secradon1} for the formal definition). \\

		Recall that one has the abelian category of
	perverse sheaves $\Perv(\PP^2)$ inside $D^b_c(\PP^2)$ and that, using perverse cohomology,  any object in $D^b_c(\PP^2)$ can be described
	by perverse sheaves (\cite{BBD}, \cite{KiehlWeissauer}).
	By Radon inversion (see Section~\ref{secradon1} and \cite{Brylinsky}, \cite{KiehlWeissauer}),
	 Brylinsky's Radon transform sends simple perverse sheaves to simple
	perverse sheaves, up to constant sheaves in $D^b_c(\tilde{\PP}^2).$ It is this result
	which makes the Radon transform useful in many applications.    \\

	 Generically, a simple perverse sheaf $K$ on $\PP^2$ is an intermediate extension $j_{!*}\cL[2]$ of an irreducible non-constant local system $\cL$ (shifted to cohomological degree $-2$)
	 along an open inclusion $j:Y\hookrightarrow \PP^2$ with  $Y$  a plane curve complement
	$Y=\PP^2\setminus C.$
		Let $\tilde{Y}=\tilde{\PP}^2\setminus \tilde{C}$ be the complement of the union $\tilde{C}$
	of the dual curve of $C$ and the lines in the dual plane, which are determined by the
	singularities of $C.$ Then the Radon transform $\cR(K)$ is locally constant
	on $\tilde{Y}$ (Prop.~\ref{prp}).  We call the local system $\cR(\cL):=\cR(K)[-2]|_{\tilde{Y}}$ the {\it Radon transform} of $\cL.$
	Under some additional assumptions which are often  fulfilled in the applications,  the Radon transform $\cR(K)$  is the  intermediate extension
	$$\cR(K)=\tilde{j}_{!*}\left(\cR(\cL)[2]\right)$$ of  $\cR(\cL)[2]$ along $\tilde{j}:\tilde{Y}\to \tilde{\PP}^2.$ Hence $\cR(\cL)$ determines $\cR(K)$ in these cases by the uniqueness of
	the intermediate extension.   	\\
		
	 The local systems
	$\cL$ and $\cR(\cL)$ are determined  in turn by their monodromy representations
	$$ \rho_\cL:\pi_1(\PP^2\setminus C)\to \GL(V)\quad \textrm{and}\quad \rho_{\cR(\cL)}:\pi_1(\tilde{\PP}^2\setminus \tilde{C})\to \GL(\tilde{V}).$$
	It is the aim of this article, to \emph{determine the monodromy representation of the
	Radon transform ${\cR(\cL)}$ in an algorithmic way} from the fundamental
	 data (= monodromy tuple of $\cL$ together with the  braid monodromy of $C$), used to describe the local system $\cL,$ see Thm.~\ref{etalem}. \\
	
		The main idea is to interprete the Radon transform
		$\cR(\cL)$ as a variation of parabolic cohomology groups over $\tilde{Y}$ in the sense of \cite{DW}
		and to use the algorithmic approach for such variations which is introduced in \cite{DW}.
		
		 We remark that
		the \'etale version of our results on the Radon transform
		 leads to new instances of the regular
		inverse Galois problem, details will appear elsewhere.

	\section{Definition and first properties of the Radon transformation}\label{secradon1}
	In the following we recall the definition of the Radon transform on
	the projective plane, following Brylinsky \cite{Brylinsky}. For results on the \'etale
	version of the
	Radon transform we refer to \cite{KiehlWeissauer}. We remark that most of the stated results also hold for higher dimensional projective spaces over more general bases
	(cf.~loc.cit.).\\
	
	{\it We  consider projective spaces and varieties  over the complex numbers. We assume additionally in Section~\ref{secradon1} that  the coefficient ring is a field
	 $k$ of characteristic zero. In the other sections, the coefficient ring   is allowed to be any field $k,$ unless otherwise stated.}\\
	
	Let $\PP^2$ be the usual projective plane and let
	$\tilde{\PP}^2$ be the dual projective plane which  parametrizes the  lines in $\PP^2.$
	Then $\tilde{\PP}^2$ is also a projective plane
	and the line in $\PP^2$ corresponding to
	$[\tilde{x}_0,\tilde{x}_1,\tilde{x}_2]\in \tilde{\PP}^2$ has the defining equation
	$$ x_0\cdot \tilde{x}_0+x_1\cdot \tilde{x}_1+x_2\cdot \tilde{x}_2=0,$$ also
	proving that $\PP^2$ is the dual of $\tilde{\PP}^2.$
	
	\begin{flushleft} {\bf Notation:}  a point $\tilde{x}\in
	\tilde{\PP}^2$ will correspond to a line  $L_{\tilde{ x}}\subseteq \PP^2.$
			By duality, an element $x\in $ $ \PP^2$ can be identified with a line
	$\tilde{L}_{x}\subseteq \tilde{\PP}^2$ in the dual plane. By construction, the line $\tilde{L}_x$  parametrizes
	the pencil of lines ${\bf L}_x$  in $\PP^2$ containing
	the point   $x.$
	\end{flushleft}
	  Let
		$$H:=\{({x},\tilde{{x}})\in \PP^2\times \tilde{\PP}^2 \mid {x}\in L_{\tilde{{x}}}\}\subset \PP^2\times \tilde{\PP}^2$$
		denote the incidence relation which is equipped  with the obvious projections $p_1:H\to \PP^2$ and
		$p_2:H\to \tilde{\PP}^2.$ \\
		
		Then one has the \emph{Radon transformation}
		$$ \cR: D^b_c(\PP^2)\to D^b_c(\tilde{\PP}^2),\, K\mapsto Rp_{2*}p_1^*(K)[1].$$
	The same construction applies to the dual space
		leading to the \emph{dual Radon transform}
		$$ \tilde{\cR}: D^b_c(\tilde{\PP}^2) \to D^b_c(\PP^2),\, \tilde{K}\mapsto Rp_{1*}p_{2}^*(\tilde{K})[1]\,.$$

	For an object $K\in D^b_c(\PP^2),$ let
$ \cR^n(K):={}^pH^n\cR(K)$ denote the $n$-th perverse cohomology of $\cR(K)$ and let $\cK:=R p_*(K),$ where $p$ denotes the structure morphism to the spectrum of $\CC.$

\begin{rem}{		The definition of the Radon transform makes sense also for more
general coefficient rings and also for more general fields of definition (cf.\cite{Brylinsky}, \cite{Patel}).
Then results whose proof depends on the decomposition theorem
(see \cite{BBD}), like the
proof of Radon inversion in \cite{Brylinsky}, \cite{KiehlWeissauer}, become more involved. This is the reason, why we assume that the coefficient ring is a field of
characteristic~$0$ in this section. We remark
that by  \cite{Patel},~Thm.~1.1 and  Rem. 1.2(3), Radon inversion still holds in our context for
coefficient rings like $\ZZ/\ell^n\ZZ$ and algebraically closed fields of definition of arbitrary characteristic.}
\end{rem}

	Let $\Perv^\mathrm{smooth}(\PP^2)$ denote the category of smooth(=constant)  perverse sheaves on $\PP^2$ and let
	$M(\PP^2)$ be the Verdier quotient $\Perv(\PP^2)/\Perv^\mathrm{smooth}(\PP^2).$ Let $M(\tilde{\PP}^2)$ be analogously defined.
Within these categories
	 one has \emph{Radon inversion} (cf.~\cite{Brylinsky}, Thm.~5.5, \cite{KiehlWeissauer},
Lemma~1.4):
	\begin{equation}\label{eqradinv} \tilde{\cR}\circ \cR\equiv \id.\end{equation}
 Especially, simple objects in $M(\PP^2)$ correspond to simple objects
	in $M(\tilde{\PP}^2)$ under Radon transform. Note that any object $K$ in $M(\PP^2)$ is represented uniquely by a \emph{reduced}
	perverse sheaf (i.e., a perverse sheaf which does not admit any nontrivial constant subobject/quotient),
	see \cite{KiehlWeissauer}, convention on Page~211. Sometimes it is useful to work directly within the categories $\Perv(\PP^2)$ and $\Perv(\tilde{\PP}^2)$ and not in their associated quotients.  This leads to the following definition:
	
	\begin{defn}\label{catP} An object $K$ in $\Perv(\PP^2)$ has the \emph{property $P$} if it has no nontrivial cohomology outside degree $0,$ i.e.,
	$\cK^i=0$ for $i\neq 0,$ where $\cK$  is as above.
	Let $\cP$ denote the category of semisimple perverse sheaves having the property $P.$ Let $\tilde{\cP}$ be the subcategory of
	$\Perv(\tilde{\PP}^2)$ defined analogously.
	\end{defn}	
\begin{rem}\label{remp} The property $P$ appears often as follows:
{\it let $K=j_{!*}(\cL)[2]$ be  an intermediate extension of an irreducible  local system $\cL$ on a dense open subset $U$  of $\PP^2$
such that $\cL$ has local monodromy  with a trivial $1$-eigenspace along a projective
line $L\subset \PP^2\setminus U,$ then $K$ is in $\cP.$}

To see this we argue as follows:
 the inclusion $j:U\to \PP^2$ is the composition of
  $j':U\to \AA^2=\PP^2\setminus L$ with $j^{''}:\AA^2\to \PP^2.$    The long exact sequence  of cohomology with compact supports with respect to the inclusion  $j^{''}:\AA^2\hookrightarrow \PP^2$ applied to the constructible
  sheaf $j_{!*}(\cL)=K[-2]$ reads
   $$ \cdots \to H^1_c(\AA^2,j'_{!*}(\cL))\to H^1(\PP^2,K[-2])\to H^1(L,K[-2]|_{L})\to H^2_c(\AA^2,j'_{!*}(\cL))\to \cdots   \;.$$
   By Artin's vanishing theorem one has  $H^1_c(\AA^2,j'_{!*}(\cL))=0,$ see \cite{Dimca}, Thm. 4.1.26 and
   Cor.~5.2.19. Since $K[-2]|_{L}$ has at most punctual support by the assumption on
   the local monodromy at $L,$ the cohomology group
   $H^1(L,K[-2]|_{L})$ also vanishes, implying that $H^1(\PP^2,K[-2])$ and hence
   $H^{-1}(\PP^2,K)=\cK^{-1}$ is zero. The other cohomology groups $\cK^{i},\, i\neq 0,$ also vanish by  the irreducibility assumption, resp.  by duality.

   We refer to
   \cite{Dimca}, Chapter 6, for other vanishing theorems in this direction. \end{rem}

\begin{prop}\label{propP} Let $K$ be an object in $\Perv(\PP^2).$ Then the following holds:
\begin{enumerate}
\item The Radon transform $\cR(K)$ is perverse if  $K$ is  in
$\cP$ or if $K$ is  a simple and non-constant perverse sheaf situated in degree $-2.$
\item
If $K$ is a simple object in $\cP,$ then $\cR(K) $ is simple as an object in $\Perv(\tilde{\PP}^2).$
\end{enumerate}
\end{prop}

\proof By \cite{KiehlWeissauer},  Lemma IV.2.2, (cf. \cite{Brylinsky}, Thm. 5.5) one has
$$ \cR^n(K)\simeq \tilde{p}^*[2]\cK^{n-1}$$ for $n<0$ and
$$ \cR^n(K)\simeq \tilde{p}^*[2]\cK^{n+1}$$ for $n>0.$ Hence the perverse sheaves  $ \cR^n(K)$ vanish for $n\neq 0,$ using that  $\cK^i=0$ for
$i\neq 0,$ proving the first claim of $(1).$ If   $K$ is  a simple and non-constant perverse sheaf situated in degree $-2,$ then it is the intermediate
extension of a non-constant irreducible local system on a dense open subset of $\PP^2.$ This implies that   $\cK^{-2}$ is trivial
(and hence also $\cK^{2}$ by duality), implying that also in this case, one has $\cR(K)=\cR^0(K)$ implying that  $\cR(K)$ is perverse. It has to be non-constant because of Radon inversion.

Let  $K$ be a simple object in $\cP.$  Since the morphism $p_2:H\to \tilde{\PP}^2$ is projective, the Radon transform $\cR(K)$ is
semisimple. The Lefschetz sequence given on page 211 of \cite{KiehlWeissauer} implies that $\tilde{p}^*[2]\cK^{-1}$ is isomorphic to the
maximal constant perverse subsheaf of $\cR^0(K)=\cR(K)$ (the latter equality follows as in the proof of $(1)$ using $\cK^{-2}=\cK^{2}=0$). Since $\cK^{-1}=0$
by property $P$ and since $\cR({K})$ is semisimple,   the Radon transform $\cR(K)$ is reduced. It follows now from Radon inversion in the form
stated in \eqref{eqradinv} that  $\cR(K)$ is again simple.
%
%
%
%
\Endproof

%
%

\section{Encoding local systems on plane curve complements}\label{secfd}
Let $C\subseteq \PP^2$ be a reduced plane algebraic curve of degree $r.$
In this section we collect the data which describe a  local system $\cL$ over an arbitrary field $k$ of rank $n\in \NN$ on
$\PP^2\setminus C,$  suitable for the Radon transformation. In the following, we will summarize the main results of
\cite{DW}, adapted to our situation:\\

 We choose a generic point ${p}_0\in\PP^2\setminus C.$  Using a projective coordinate change on $\PP^2$ we
can assume that $p_0=[1:0:0].$  Using the conventions from the beginning of Section~\ref{secradon1},
  we write  $\tilde{L}_{{p}_0}\subseteq \tilde{\PP}^2$ for the corresponding  line in the dual plane.  Hence $$\tilde{L}_{p_0}=\{[0:\tilde{x}_1:\tilde{x}_2]\mid [\tilde{x}_1:\tilde{x}_2]\in \PP^1\}.$$
%
  We also  choose an element $\tilde{p}_0\in \tilde{L}_{p_0}\subseteq \tilde{\PP}^2$ such that $L_{\tilde{p}_0}\subseteq {\PP}^2$ intersects $C$ transversally, where $L_{\tilde{p}_0}$ is the line in $\PP^2$ determined by $\tilde{p}_0. $
   By possibly using a further coordinate change fixing $p_0$ we can assume that
  $\tilde{p}_0=[0:1:0]\in \tilde{\PP}^2.$ \\
%
%

  As in the previous section, the incidence relation is defined as
  $$H:=\{([x_0:x_1:x_2],[\tilde{x}_0:\tilde{x}_1:\tilde{x}_2])\in \PP^2\times \tilde{\PP}^2 \mid
 x_0\cdot \tilde{x}_0+x_1\cdot \tilde{x}_1+x_2\cdot \tilde{x}_2=0 \}\subset \PP^2\times \tilde{\PP}^2,$$
  coming
		with its projections  $p_1:H\to \PP^2$ and
		$p_2:H\to \tilde{\PP}^2.$
  We can view the pencil  ${\bf L}_{p_0}$ as a subvariety of $H:$
  $${\bf L}_{p_0}=p_2^{-1}(\tilde{L}_{p_0})=\{([x_0:x_1:x_2],[0:\tilde{x}_1:\tilde{x}_2])\in H\}   .$$

     Let ${\tilde{y}_1},\ldots,{\tilde{y}_s}\in \tilde{L}_{p_0}$ denote the points which are \emph{exceptional } with respect to
$C,$ meaning that each $L_{\tilde{ y}_i}$ is either tangent to  $C$ or that $L_{\tilde{ y}_i}$ contains a singularity of $C.$
By possibly using another projective coordinate change (which fixes $p_0$ and fixes $L_{\tilde{p}_0}$ setwise), we can assume that $\tilde{y}_s=[0:0:1].$

One has an isomorphism of $\PP^1$-bundles
\begin{equation}\label{eqLambda}  \Lambda:{\bf L}_{p_0}\setminus p_2^{-1}({\tilde{y}_s})\to \PP^1_{\AA^1}=\PP^1\times \AA^1\end{equation} given by
  $$ ([x_0:x_1:x_2],[0:1:\tilde{x}_2]) \to ([x_0:x_1],\tilde{x}_2),$$ where the bundle maps are the respective second projections.

  Let $$S:=\tilde{L}_{p_0}\setminus \{{\tilde{y}_1},\ldots,{\tilde{y}_s}\},$$ let $X:=p_2^{-1}(S)$
  and $X_0:=p_2^{-1}(\tilde{p}_0).$ Then
 $D:=p_1^{-1}(C) \cap X$ is a divisor of relative degree $r$ over $S$ (denoted $r$-configuration in \cite{DW}, Section~2.1).
 Let $U:=X\setminus D,$ write $\pi:U\to S$ for the restriction of $p_2$ to $U,$ and let $U_0:=p_2^{-1}({\tilde{p}_0})\cap U.$
 Then $\pi:U\to S$ is a topological fibration with typical fibre $U_0$ and hence yields a short exact sequence as follows (cf.~\cite{DW}, Formula~(6)):
 \begin{equation}\label{dw1}
 1\to \pi_1(U_0,(p_0,\tilde{p}_0))\to \pi_1(U,(p_0,\tilde{p}_0))\to \pi_1(S,\tilde{p}_0)\to 1\,.
 \end{equation}
 By restricting the above $\PP^1$-bundle  isomorphism
 $\Lambda: {\bf L}_{p_0}\setminus p_2^{-1}({\tilde{y}_s})\to \PP^1_{\AA^1}$ to $X$ we obtain
 a trivializing  isomorphism $\lambda: X\simeq \PP^1_S$ of $\PP^1$-bundles
 which identifies the map $p_2:X\to S$ with
 the structural map $\PP^1_S\to S.$ One has $$\lambda((p_0,[0,1,\tilde{x}_2]))=(\infty, \tilde{x}_2),$$ where $\infty=[1:0]$ is the point at infinity in $\PP^1.$
  Moreover,   since $p_0\notin C$ as $p_0$ is generic w.r. to $C,$    the section $\infty \times S=\lambda(p_1^{-1}(p_0))$ is disjoint from
 $\lambda(D).$   \\

  As in \cite{DW}, Section 1.2, we fix a homeomorphism (called a \emph{marking}) $\kappa: X_0\to \PP^1(\CC),$ mapping
  $P_0:=X_0\setminus  D$ to a subset $\{x_1,\ldots,x_r\}$ of the real line inside the Riemann sphere $\PP^1(\CC)$
  and $p_0$ to an element $x_0$ in the upper half plane. Let us, just for the moment, identify
  $U_0=X_0\setminus P_0$ with $\PP^1(\CC)\setminus \{x_1,\ldots,x_r\}$ via $\kappa.$
  As  in \cite{DW} now  choose counterclockwise simple loops $\gamma_1,\ldots,\gamma_r$ around $x_1,\ldots,x_r$ which are  based at $x_0$ as follows: Each $\gamma_i,\,
  i=1,\ldots,r$ is a simple closed loop which intersects the real line exactly twice, first the open interval $(x_{i-1}, x_i),$  then the interval $(x_i, x_{i+1})$ (where for $i=1,$ the interval
  $(x_{i-1},x_i)$ is iterpreted as $(-\infty,x_1)$ and for $i=r$ the interval $(x_{i},x_{i+1})$
  is interpreted as $(x_{i+1},+\infty)$).

  By taking inverse images under $\kappa,$ this  leads to a  standard presentation
 of the  fundamental group
   $$\pi_1(U_0,(p_0,\tilde{p}_0))=\langle \gamma_1,\ldots,\gamma_r\mid \gamma_1\cdots \gamma_r=1\rangle.$$  We remark that (up to a M\"obius transform) the construction
   of the $\gamma_1,\ldots,\gamma_r$ is the same as the construction
   of a geometric basis of meridians
   in \cite{Bartolo}. \\

  By the above presentation of $\pi_1(U_0),$ the restriction $\cL_0:=\cL|_{U_0}$   is  determined   by the tuple of matrices
  $$\g_\cL:=(g_1,\ldots,g_r)\in \GL(V)^r,\quad g_1\cdots g_r=1,$$ where $\gamma_i$ is sent to $g_i$ under the monodromy representation $\rho_{\cL}$ of $\cL$
  (here $V\simeq k^n$). We call $\g$ the \emph{monodromy tuple} of $\cL_0$ (resp. of $\cL$) with respect to $ \gamma_1,\ldots,\gamma_r.$

\begin{rem} \label{dwrem}
By Lefschetz' hyperplane theorem (see \cite{GoreskyMacPherson}, \cite{dimca}), the fundamental group
$\pi_1(U_0)$ maps surjectively onto $\pi_1(\PP^2\setminus C).$ Hence the monodromy tuple $\g$ of $\cL_0$ also
determines the local system $\cL$ up to isomorphism. \end{rem}

  Let $$\cO_r:=\{P\subseteq \AA^1\mid \#P=r\}$$ denote the configuration space parametrizing subsets of $\CC$ of cardinality $r$ and let
 $$\cO_{r,1}:=\{(q,P)\in \PP^1\times \cO_r\mid q\notin P\}.$$

 As explained in \cite{DW},  Section 2.2, our choice of a marking
 defines generators $\beta_1,\ldots,\beta_{r-1}$ of the braid group
 $ \cB_r:=\pi_1(\cO_r,P_0)$ as follows: under $\kappa,$ a representative of $\beta_i$ fixes all points different from $\{x_i, x_{i+1}\}$
 and the image under $\kappa$ of $x_i$ moves to $x_{i+1}$ in the  lower half plane
  and  (the image under $\kappa$ of) $x_{i+1}$ travels to $x_i$ in the upper half plane.
This leads to the standard  presentation of $ \cB_r:$  $$ \langle \beta_1,\ldots,\beta_{r-1}\mid \beta_{i}\beta_{i+1}\beta_{i}=  \beta_{i+1}\beta_{i}\beta_{i+1} \textrm{ for } i=1,\ldots,r-1 \textrm{ and }  \beta_i\beta_j=\beta_j\beta_i \textrm{ for }
	|i-j|>1\rangle.$$
	
	Let $\cB_{r,1}:=\pi_1(\cO_{r,1},(\infty,P_0)).$
	The natural projection $\cO_{r,1}\to\cO_r$ is a topological fibration with typical fibre
 $U_0$, and admits a section $P\mapsto(\infty,P).$  It yields a
split exact sequence of fundamental groups
\begin{equation} \label{splitseq1}
   1 \;\to\; \pi_1(U_0,{\infty}) \;\To\; \pi_1(\cO_{r,1},(\infty,P_0))
     \;\To\; \cB_r \;\to\; 1.
\end{equation}
We may identify $\cB_r$ with its image in $\pi_1(\cO_{r,1})$ under the
splitting induced from the above section.  Then the conjugation operation of
$\cB_r$  on $\pi_1(\PP^1\setminus P_0,{\infty})$ is as follows (\cite{DW}, Formula~10):
	\begin{equation} \label{locals5eq2}
  \beta_i^{-1}\,\gamma_j\beta_i \;=\quad
  \begin{cases}
    \quad \gamma_i\gamma_{i+1}\gamma_i^{-1}, &
                              \quad\text{\rm for $j=i$,} \\
    \quad \gamma_i, & \quad \text{\rm for $j=i+1$,} \\
    \quad \gamma_j, & \quad \text{\rm otherwise.}
  \end{cases}
\end{equation}


 Since $\lambda(D)$  is disjoint from $\infty \times S$ we obtain  a map $p:S\to \cO_r$
 by  sending  $s$ to the fibre of $\lambda(D)$ over $s.$ Let $\varphi:\pi_1(S,\tilde{p}_0)\to \cB_r$ be the homomorphism induced by $p.$ Then the above short exact sequence \eqref{dw1} can be viewed as the pullback of the short exact sequence
 \eqref{splitseq1} along $\varphi$ (cf.~\cite{DW}, Section 2.3, for these results). Using the splitting of \eqref{dw1} coming from the
 $\infty$-section, one can view $\pi_1(S)$ as a subgroup of $\pi_1(U).$ By construction, the resulting conjugation action
 of $\pi(S)$ on $\pi_1(U_0)$ factors through the map $\varphi$  and is determined by Formula~\eqref{locals5eq2}. \\

		\begin{defn} Let $\tilde{\gamma}_1,\ldots,\tilde{\gamma}_s$ be simple loops in $S$
				satisfying the product
		rule $\tilde{\gamma}_1\cdots \tilde{\gamma}_s=1,$ constructed using a marking as above.
				Then their images under $\varphi$ are words in the braids $\beta_1,\ldots, \beta_{r-1}:$
		$$\varphi(\tilde{\gamma}_i)=\omega_i(\beta_1,\ldots,\beta_{r-1})=\omega_i\in \cB_r.$$
		These words constitute the \emph{braid monodromy} of $C.$
		\end{defn}

	 \begin{defn} The \emph{fundamental data} of a local system $\cL$ on $\PP^2\setminus C$ consist
	 of the following:
	 	 \begin{itemize}
	 \item The monodromy tuple  $(g_1,\ldots,g_r)\in \GL(V)^r\,(V\simeq k^n)$ of $\cL$ with respect to
	 fixed generators
	 $\gamma_1,\ldots,\gamma_r$ of $\pi_1(U_0).$
	 \item The braid monodromy $(\omega_1,\ldots, \omega_s)\in \cB_r^s$
	 of $C$     with respect to fixed generators
	 $\tilde{\gamma}_1,\ldots,\tilde{\gamma}_s$ of $\pi_1(S).$
	 \end{itemize}
	 \end{defn}

			\begin{rem} \label{remvankampen} \begin{enumerate}
		\item
		As mentioned above, the trivializing morphism $\lambda: X\to \PP^1_{S}$ has
		the property that   $\lambda(D)$ and $ \infty\times S$ are disjoint. It follows that the restriction
		of $\lambda$ to $U\setminus \lambda^{-1}(\infty\times S)$ constitutes an isomorphism
		$$U\setminus \lambda^{-1}(\infty\times S) \to \AA^2\setminus (C\cup L_{\tilde{y}_1}\cup \ldots \cup
		L_{\tilde{y}_{s-1}})$$ under which the map
		$U\setminus \lambda^{-1}(\infty\times S)\to S$
		is transformed into the second projection
		$$ \pr_2:\AA^2\setminus (C\cup L_{\tilde{y}_1}\cup \ldots \cup
		L_{\tilde{y}_{s-1}})\to \AA^1_y\setminus \{-{y_1},\ldots,-{y_{s-1}}\},\, (x,y)\mapsto y,$$
		composed by the map
		$$ -{\bf 1}:\AA^1_y\setminus \{-{y_1},\ldots,-{y_{s-1}}\} \to \AA^1_y\setminus \{{y_1},\ldots,{y_{s-1}}\} ,\,
		y\mapsto -y.$$ Hence the braid monodromy in the above sense coincides, up to a twist by the morphism $-
		{\bf 1},$
		indeed with the classical notion of
		braid monodromy (cf.~\cite{dimca}).
		\item \label{item3}  In practise, we
		first choose generators $-\tilde{\gamma}_1,\ldots,-\tilde{\gamma}_s$ of  $\pi_1(\AA^1_y\setminus \{-{y_1},\ldots,-{y_{s-1}}\})$
		satisfying the product
		rule $(-\tilde{\gamma}_1)\cdots (-\tilde{\gamma}_s)=1$ and define
		$\tilde{\gamma}_1,\ldots,\tilde{\gamma}_s$ to be their images under $\lambda^{-1}\circ -{\bf 1}.$ Then the classical braid monodromy of the plane curve complement yields the braid monodromy defined above.
		\item If $\tilde{\gamma}_1,\ldots,\tilde{\gamma}_s$ are constructed as in (\ref{item3}),
		by van Kampens theorem, the fundamental group $\pi_1(\PP^2\setminus C)$
			 has a presentation as follows (where the conjugation operation of $\omega_i$ on $\gamma_i$ is as in
			 \eqref{locals5eq2}):
			$$  \langle \gamma_1,\ldots,\gamma_r\mid \gamma_1\cdots\gamma_r=1,\,
			{\omega_i}^{-1}\gamma_j\omega_i=\gamma_j \textrm{ for } i\in \{1,\ldots,s\}, j\in \{1,\ldots,r\}\rangle.$$ Hence the monodromy tuple $(g_1,\ldots,g_r)\in \GL(V)^r$ of $\cL$ also satisfies these relations.

					\end{enumerate}
			\end{rem}

%

\section{Deformations of local systems and cocycles
using the braid group}

		We continue by summarizing the results of \cite{DW} concerned with deformations
		of local system along braids: Using the notions of the previous section, define
\[
    \E_r \;:=\; \{\;\g=(g_1,\ldots,g_r) \mid
                  g_i\in\GL(V),\;\; \prod_ig_i=1 \;\}.
\]
As described in the previous Section, an
element $\g\in\E_r$ corresponds to a representation
$\rho_0:\pi_1(U_0)\to\GL(V)$ (set $\rho_0(\gamma_i):=g_i$) and
hence to a local system  on $U_0$.  \\

As described in \cite{DW}, Section 2.2, the natural deformation operation of the braid group
on local systems (by composing $\rho_0$ with the operation of  $\beta \in \cB_r$ on $\pi_1(U_0)$)
amounts to a right-action of $\cB_r$ on  $\E_r$ as follows: \begin{equation}  \label{locals5eq4}
   \g^{\beta_i} \;=\; (g_1,\ldots,g_{i+1},g_{i+1}^{-1}g_ig_{i+1},\ldots,g_r).
\end{equation}

This deformation operation enables the computation of the braiding operation
on cocycles which will be considered in the next section:
given $\g\in\E_r$, we have  the vector space
\[
    H_\g \;=\; \{\;(v_1,\ldots,v_r) \mid v_i\in{\rm Im}(g_i-1),\;
       v_1\cdot g_2\cdots g_r+v_2\cdot g_3\cdots g_r+\cdots +v_r  =0   \;\}\subset V^r.
\]
We remark that we have a natural  isomorphism between $H_\g/E_\g$ (where $E_\g\leq H_\g$ is as defined below) and the
parabolic cohomology group $H^1(\PP^1,j_*\cL)$ (where $j$ denotes the inclusion of
$U_0$ to $\PP^1,$ cf. \cite{DW}, Section~2).

In loc.cit., Section~3, for $\g\in \E_r$ and $\beta\in \cB_r$ a linear map \[
      \Phi(\g,\beta):\,H_{\g} \;\To\; H_{\g^{\beta}}.
\]
is defined as follows (recall our convention of linear maps acting from the right): for a standard generator
$\beta_i\in \cB_r$ one has
\begin{multline} \label{locals5eq6}
 \qquad (v_1,\ldots,v_r){\Phi(\g,\beta_i)} \\
    \;=\;  (v_1,\ldots,v_{i+1},\,
    \underbrace{v_{i+1}(1-g_{i+1}^{-1}g_ig_{i+1})+v_ig_{i+1}}_{
        \text{\rm $(i+1)$th entry}},\,\ldots,v_r).
\end{multline}
Note that both vector spaces $H_{\g} $ and $H_{\g^{\beta}}$ are subspaces of $V^r.$ Hence
${\Phi(\g,\beta_i)}$ is  given by restriction of the following endomorphism $M(\g,\beta_i)\in \End(V^r)$
in matrix form to $H_{\g} :$

\begin{equation}\label{Phi}
M(\g,\beta_i)=\left(\begin{array}{llllllllll}
1&&&\\
&&\ddots&&&\\
&&&\quad 1&&\\
&&&&\,\,\,0&\quad g_{i+1}&\\
&&&&\,\,\,1&\quad (1-g_{i+1}^{-1}g_ig_{i+1})&&\\
&&&&&&&1&&\\
&&&&&&&&\quad \ddots&\\
&&&&&&&&&\quad \quad 1
\end{array}\right)\,.
\end{equation}
Moreover, one has the `cocycle rule' which allows one to compute $\Phi(\g,\beta)$ for a general element $\beta\in \cB_r$ in terms of the $\Phi(\g,\beta_i):$

\begin{equation} \label{locals5eq7}
    \Phi(\g,\beta'\beta'')  \;=\; \Phi(\g,\beta')\cdot\Phi(\g^{\beta'},\beta'') \quad
    (\beta',\beta''\in \cB_r ).
\end{equation}
The product on the right hand side of \eqref{locals5eq7} is defined as
the linear map from $H_{\g}$ to $H_{\g^{\beta'\beta''}},$ acting from the right,  obtained by
first applying $\Phi(\g,\beta')$ and then
$\Phi(\g^{\beta'},\beta'')$. \\

The cocycle rule especially implies that
$${\Phi(\g,\beta_i^{-1})}=\Phi(\g^{\beta_i^{-1}}, \beta_i)^{-1}$$ is given by restriction of
 $$M(\g^{\beta_i^{-1}},\beta_i)^{-1}\in \End(V^r)$$  to $H_{\g} .$ The latter endomorphism
 is denoted by $M(\g,\beta_i^{-1})$ for convenience. By  the definition
 of  $M(\g^{\beta_i^{-1}},\beta_i)$ and inversion one derives the following matrix-expression:
  \begin{equation}\label{Phi2}
M(\g,\beta_i^{-1})=\left(\begin{array}{llllllllll}
1&&&\\
&&\ddots&&&\\
&&&\quad 1&&\\
&&&&\,\,\,(g_{i+1}-1)g_i^{-1}&\quad 1&\\
&&&&\,\,\,\quad\quad  g_i^{-1}&\quad 0&&\\
&&&&&&&1&&\\
&&&&&&&&\quad \ddots&\\
&&&&&&&&&\quad \quad 1
\end{array}\right)\,.
\end{equation}

\begin{prop}\label{rembraids}\begin{enumerate}
\item
Let $\omega= \beta_{i_1}^{\epsilon_1}\cdots \beta_{i_k}^{\epsilon_k}\in \cB_r$ be a word in the standard generators $\beta_1,\ldots,\beta_{r-1}$
of $\cB_r$ (with $\epsilon_1,\ldots,\epsilon_k\in \{\pm 1\}$) and let $\g\in \E_r.$ Then the isomorphism ${\Phi(\g,\omega)}:H_{\g}\to H_{\g^\omega}$ is given by the restriction of the following endomorphism  $M(\g,\omega)\in \End(V^r)$ to
$H_{\g}:$
$$M(\g,\omega)=M(\g,\beta_{i_1}^{\epsilon_1})\cdot M(\g^{\beta_{i_1}^{\epsilon_1}},\beta_{i_2}^{\epsilon_2})\cdots M(\g^{\beta_{i_1}^{\epsilon_1}\cdots \beta_{i_{k-1}}^{\epsilon_{k-1}}},\beta_{i_k}^{\epsilon_k}).$$
%
\item
Let
\[
      E_{\g} \;:=\; \{\;(v\cdot(g_1-1),\ldots,v\cdot(g_r-1)) \;\mid\;
                           v\in V \;\}\,.
\]  Then $\Phi(\g,\omega)$  induces an isomorphism
\[
       \bar{\Phi}(\g,\omega):\,W_{\g}:=H_{\g}/E_{\g} \;\liso\;
              W_{\g^{\omega}}.
\]	
\end{enumerate}
\end{prop}

\proof The first claim follows from the cocycle rule \eqref{locals5eq7}. The second claim follows from the fact that
$\Phi(\g,\omega)$ maps $E_\g$ to $E_{\g^\omega},$ cf.~\cite{DW}, Section~2.2.\Endproof

		\section{The monodromy of the Radon transform}	\label{Secradtrans}

					 We assume in the following, {\it that $K=j_{!*}\cL[2]$ with $\cL$ an irreducible and non-constant
	local system of rank $n$ 	on $Y=\PP^2\setminus C,$ where $C$ is a reduced plane curve of degree $r.$}\\

Using a suitable coordinate change, we can assume that we are  in a situation as described in Section~\ref{secfd}.
We will freely use the notation of the previous sections. Especially, let $$\left( (g_1,\ldots,g_r), (\omega_1,\ldots, \omega_s)\right)\in \GL(V)^r\times \cB_r^s$$
	be the fundamental data of $\cL$ with respect to a fixed choice of generators
	 $$\gamma_1,\ldots, \gamma_r\in \pi_1(U_0,(p_0,\tilde{p}_0))\,\,\textrm{ and } \,\,\tilde{\gamma}_1,\ldots,\tilde{\gamma}_s\in \pi_1(S,\tilde{p}_0)$$ as in the previous sections. \\

		Let $\tilde{C}\subseteq \tilde{\PP}^2$ be the union of the dual curve of $C$ and the lines in
		$\tilde{\PP}^2$ which correspond to pencils of lines in $\PP^2,$ passing through the singularities of $C.$ Moreover, let
		$\tilde{Y}:=\tilde{\PP}^2\setminus \tilde{C}.$
		For $\tilde{x}\in \tilde{\PP}^2$ let $j_{\tilde{x}}: L_{\tilde{x}}\setminus C\hookrightarrow L_{\tilde{x}}$ be  the natural inclusion and let $\cL_{\tilde{x}}:=\cL|_{L_{\tilde{x}}\setminus C}.$

		\begin{prop}\label{prp} For $\tilde{x}\in \tilde{Y}$ let $ \cR(K)_{\tilde{x}}$ be the stalk of the Radon transform at $\tilde{x}$ in the derived sense.  Consider the \emph{parabolic cohomology group}
		$$ H^1_{par}(\cL_{\tilde{x}}):=H^1(L_{\tilde{x}}, j_{\tilde{x}*}(\cL_{\tilde{x}})).$$ Then there is an isomorphism
		$$ \cR(K)_{\tilde{x}}[-2]\simeq H^1_{par}(\cL_{\tilde{x}}).$$ Moreover, the dimension of $H^1_{par}(\cL_{\tilde{x}})$ is
			\begin{equation}\label{degree}     n(r-2)-\sum_{i=1}^r \dim(V^{g_i}).\end{equation}
		Especially, 	the restriction of $\cR(K)$ to $\tilde{Y}$ is
	a local system whose rank is given by \eqref{degree},  placed in cohomological degree $-2.$
\end{prop}
		
	\proof 	If $\tilde{x}\in \tilde{Y},$ then  $\pi_1(L_{\tilde{x}}\setminus C)$ maps surjectively onto  $\pi_1(U)$ (see~\cite{dimca}), implying that
	 $\cL_{\tilde{x}}=\cL|_{L_{\tilde{x}}\setminus C}$
		is an irreducible and  non-constant local system on $L_{\tilde{x}}\setminus C.$ Hence $H^i(L_{\tilde{x}}, j_{\tilde{x}*}(\cL_{\tilde{x}}))=0$ for $i\neq 1.$
				By proper basechange and the irreducibility assumption, the stalk $\cR(K)_{\tilde{x}}$
		of the Radon transform $$\cR(K)=Rp_{2*}p_1^*(j_{!*}\cL[2])$$  is the cohomology group $H^1(p_2^{-1}(\tilde{x}),p_1^{*}(j_{!*}(\cL))|_{p_2^{-1}(\tilde{x})}),$ placed in degree $-2.$
		The first projection induces an isomorphism $p_1: p_2^{-1}(\tilde{x})\to L_{\tilde{x}}$ with $$p_1^*( j_{\tilde{x}*}(\cL_{\tilde{x}}))\simeq \left(p_1^{*}(j_{!*}(\cL))\right)|_{p_2^{-1}(\tilde{x})}.$$ Hence
		$$\cR(K)_{\tilde{x}}[-2]\simeq H^1(L_{\tilde{x}}, j_{\tilde{x}*}(\cL_{\tilde{x}})). 	$$
		It is a consequence of  the Ogg-Shafarevich Formula (cf.~\cite{DW}), that for all $\tilde{x}\in \tilde{Y}$ the dimension of  $H^1_{par}(\cL_{\tilde{x}})$ is as in Formula~\eqref{degree}, since the local monodromy is constant along the components of $C$ and since each line $L_{\tilde{x}}$ intersects $C$ transversally.
	 This implies that the restriction of $\cR(K)$ to $\tilde{Y}$ is
	a local system whose degree is given by Eq.~\eqref{degree},  placed in cohomological degree $-2.$  \Endproof
			
		\begin{defn} \label{def1} 	Let  $K$ and $\cL$ be as above.	The \emph{Radon transformation} of $\cL$ is the local system
		$$\cR(\cL):= (\cR(K)|_{\tilde{Y}})[-2]$$ on $\tilde{Y}$.
		\end{defn}

\begin{rem}	It follows from  Prop.~\ref{prp} that $\cR(\cL)$ is the zero-local system if
	$$   n(r-2)-\sum_{i=1}^r \dim(V^{g_i})=0.$$  However, if
	$$   n(r-2)-\sum_{i=1}^r \dim(V^{g_i})\neq 0,$$ then usually
	$\cR(\cL)$ is non-constant and hence determines $\cR(K)$ uniquely.
	\end{rem}

%
%
%
		
%


		An immediate consequence of Lefschetz' hyperplane theorem is the following:
	\begin{lem} Let $S=\tilde{L}_{p_0}\setminus \{\tilde{y}_1,\ldots,\tilde{y_s}\}\subseteq \tilde{Y}.$ Then the
	Radon transform  $\cR(\cL)$   is determined  by its restriction to
	$S.$
	Then $\cR(\cL)$ is determined up to isomorphism by the monodromy tuple $$
	\tilde{\g}=(\tilde{g}_1,\ldots,\tilde{g}_s)\in \GL(\cR(\cL)_{\tilde{p}_0})^s$$ of $\cR(\cL) $ with respect to $\tilde{\gamma}_1,\ldots,\tilde{\gamma}_{s}.$\Endproof
		\end{lem}

		Let $H$ be the incidence relation used in the definition of the Radon transform, equipped with its
		projections $p_1:H\to \PP^2$ and $p_2:H\to \tilde{\PP}^2.$
		Let $\cL_0=\cL_{\tilde{p}_0}$ be the restriction of $p_1^*\cL$ to $U_0\subseteq U,$ with
		$U_0\subseteq U\subseteq X\subseteq H$ as in Section~\ref{secfd}.
		Then $p_1^*\cL$ is a variation of $\cL_0$ over $\tilde{Y}$ with $\cR(\cL)$ the parabolic cohomology
		of
		this variation (in the sense of \cite{DW}, Def.~2.1).
		
		Recall from \cite{DW} that the parabolic cohomology group
		$H^1_{par}(\cL_0)=\cR(\cL)_{\tilde{p}_0}$ is isomorphic to a subspace of  the first group
		cohomology  $H^1(\pi_1(\AA^1\setminus P_0),V)\,(P_0=X_0\setminus U_0)$ where $V\simeq k^n$ is the stalk of $\cL_0$
		at $(p_0,\tilde{p}_0)\in U_0$ and $V$ is a $\pi_1(\AA^1\setminus P_0)$-module via the monodromy representation
	$\rho_{\cL_0}$	of $\cL_0.$ The group cohomology $H^1(\pi_1(\AA^1\setminus P_0),V)$ is formed, up to equivalence induced by boundaries, by maps (called cocycles) $$\delta\colon\pi_1(\AA^1\setminus P_0)\to V$$ satisfying the {\it  cocycle rule }
	\begin{equation}\label{cocyclerule}	
		\delta(\gamma\gamma')=\delta(\gamma)\rho_{\cL_0}(\gamma')+\delta(\gamma') .\end{equation}

\begin{lem}\label{lemH1} Let $\g:=(g_1,\ldots,g_r)\in \GL(V)^r$ be the monodromy tuple of $\cL_0.$
Let $W_{\g}:=H_{\g}/E_{\g}$ be as in Prop.~\ref{rembraids}(ii).
	Then the evaluation map  $$ \cR(\cL)_{\tilde{p}_0}=H^1_{par}(\cL_0)\to W_{\g}, \quad [\delta]\mapsto
	(\delta(\gamma_1),\ldots,\delta(\gamma_r))+E_{\g}$$ is an isomorphism. \end{lem}
\proof Follows from
		\cite{DW}, Section~1.3. \Endproof

%

%
%



	\begin{thm} \label{etalem} Let $\cL$ be an irreducible and non-constant local system on a plane curve complement
	$\PP^2\setminus C$ with fundamental data
	$$\left( \g=(g_1,\ldots,g_r), (\omega_1,\ldots, \omega_s)\right)\in \GL(V)^r\times \cB_r^s$$ with
	respect to a fixed choice of generators
	 $$\gamma_1,\ldots, \gamma_r\in \pi_1(U_0,p_0)\,\,\textrm{ and } \,\,\tilde{\gamma}_1,\ldots,\tilde{\gamma}_s\in \pi_1(S,\tilde{p}_0).$$
	
	Then the monodromy tuple
	 $$\tilde{\g}=(\tilde{g}_1,\ldots, \tilde{g}_s)\in \GL(W_{\g})^s$$ of the Radon transform $\cR(\cL)$ with respect to
	$\tilde{\gamma}_1,\ldots,\tilde{\gamma}_s$ is as follows:
	 For $i\in \{1,\ldots ,s\},$ write
	 $$\omega_i= \beta_{i_1}^{\epsilon_1}\cdots \beta_{i_k}^{\epsilon_k}\in \cB_r$$ as a word in the standard generators $\beta_1,\ldots,\beta_{r-1}$
of $\cB_r$ with $\epsilon_1,\ldots,\epsilon_k\in \{\pm 1\}.$
Let   $$M(\g,\omega_i)=M(\g,\beta_{i_1}^{\epsilon_1})\cdot M(\g^{\beta_{i_1}^{\epsilon_1}},\beta_{i_2}^{\epsilon_2})\cdots M(\g^{\beta_{i_1}^{\epsilon_1}\cdots \beta_{i_{k-1}}^{\epsilon_{k-1}}},\beta_{i_k}^{\epsilon_k})\in \End(V^r)$$ be as in Prop.~\ref{rembraids}.
	 Then  $\tilde{g}_i$ is induced by the operation of
	 $M(\g,\omega_i)\in \End(V^r)$ on the subquotient  $W_\g=H_\g/E_\g\leq V^r/E_\g.$
%
%
%
%
\end{thm}
	
	\proof Consider the second projection
	$p_2:H\to \tilde{\PP}^2,$ where $H$ is the incidence relation, and  let
	$S=\tilde{L}_{p_0}\setminus \{\tilde{y}_1,\ldots,\tilde{y}_s\} \subseteq \tilde{\PP}^2. $
	Let $X=p_2^{-1}(S)$ and let $D=p_1^{-1}(C)\cap X$ be as in the previous section.
	 Then the  pair $(X,D)$ is an \emph{$r$-configuration}
	in the sense of \cite{DW}, Section~2.1 (i.e., for $s\in S,$ the fibres $X_s$ are  Riemann surfaces of genus $0$
	and $\#D\cap X_s=r$). Recall that $U$ was defined to be   $X\setminus D.$
	
	There exists an isomorphism $\lambda:X\to \PP^1_S$ such
	that $\lambda(D)$ is disjoint from $\infty\times S.$
	Therefore we can embed $S$ into $U\subseteq X$ via $\lambda^{-1}(\infty\times S).$ This
	induces  a section $s:\pi_1(S)\to \pi_1(U)$
	to the homomorphism
	$\pi_1(U)\to \pi_1(S)$ from the short exact sequence in \eqref{dw1}.  Write
$\chi:\pi_1(S)\to \GL(V)$ for the composition $\rho_{	p_1^*(\cL)}\circ s.$

For $h\in \GL(V)$ and $\g^h\in \GL(V)^r$ the tuple obtained by simultaneously conjugating
the components of $\g$ by $h,$ one gets an isomorphism $$\bar{\Psi}(g,h):W_{\g}\to W_{\g^h}$$
from the diagonal operation of $h$ on $V^r.$
Since $p_0$ is a smooth point of $\cL$ and since $$p_1^{-1}(p_0)=\lambda^{-1}(\infty\times S),$$
 the homomorphism $\chi$ is trivial, implying
	that $\bar{\Psi}(g,\chi(\tilde{\gamma}_i))$ is trivial.

	As already mentioned, the divisor
	 $\lambda(D)$ is disjoint from $\infty\times S$ (hence $\lambda$ is an \emph{affine frame}
	 as defined in loc.cit., Section~2.4). 	
	Hence, by \cite{DW}, Thm.~2.5, under the identification of $W_\g=\cR(\cL)_{\tilde{p}_0},$
	 the monodromy representation $\rho_{\cR(\cL)}: \pi_1(S)\to \GL(W_\g)$ is as follows:
	 $$ \rho_{\cR(\cL)}(\tilde{\gamma}_i)=\bar{\Phi}(\g,\omega_i)\cdot \bar{\Psi}(g,\chi(\tilde{\gamma}_i))=\bar{\Phi}(\g,\omega_i)\in \GL(W_\g),$$ where
	$ \bar{\Phi}(\g,\omega_i)$ is as in Prop.~\ref{rembraids}. \Endproof

%
	
		\begin{rem}	Let $K$ be a simple and non-constant perverse sheaf whose coefficient field $k$ has characteristic $0.$ It is a consequence  of Radon inversion
		that		the Radon transform $\cR(K)$ has a unique
		simple non-constant subfactor $S$ (cf. the discussion in Section~\ref{secradon1}). 	Let us additionally
		assume that the following generic condition holds: the Radon transform  $\cR(\cL)$ is a non-constant local system on
		 $\tilde{Y}.$
Then $j_{!*}(\cR(\cL))[-2]$ contains this unique subfactor $S.$ If moreover $K$ satisfies the condition $P$
then Prop.~\ref{propP} implies that
$$ j_{!*}(\cR(\cL))[-2]=\cR(K).$$
  The local system $\cR(\cL)$ in turn is  determined up to isomorphism  by the monodromy tuple
  (with respect to $\tilde{\gamma}_1,\ldots,\tilde{\gamma}_s$), which is algorithmically determined by Thm.~\ref{etalem}.
 \end{rem}

	\section{Examples: an elliptic fibration and the classical  Zariski pair}
Let $k$ be a field of characteristic $\neq 2.$  Let $$
\tilde{x}_1:=[2:-4:3],\, \tilde{x}_2:=[1,-1,1],\, \cx_3:=[1:0:0],\, \cx_4:=[1:0:-1]\, \in \tilde{\PP}^2$$
and let $\cL$ be the rank-one local system on $$\PP^2\setminus C\quad  (C:=\bigcup_{i=1}^4 L_{\cx_i})$$ having quadratic
local
monodromy $-1\in k^*$ along $L_{\tilde{x}_i}\, (i=1,2,3,4).$
Hence the monodromy tuple of $\cL$
(defined as in Section~\ref{secfd}) is
$\g=(-1,-1,-1,-1)\in (k^*)^4.$ \\

We remark that $C$ represents a generic quadruple of projective lines in $\PP^2.$ \\

The singularities of $C$ are
$$ x_1:=[-\frac{1}{2}:\frac{1}{2}:1],\, x_2:=[0:\frac{3}{4}:1],\, x_3:=[0:1:1],$$
$$  x_4:=[1:\frac{5}{4}:1],\,
x_5:=[1:2:1],\, x_6:=[0:1:0].$$
Hence the Radon transform $\cR(\cL)$ of $\cL$ is a $k$-local system on $\tilde{\PP}^2\setminus \tilde{C},$
where $\tilde{C}$ is the union of $\tilde{L}_{x_1},\ldots,\tilde{L}_{x_6}.$

\begin{prop}\label{prop4}
The Radon transform $\cR(\cL)$ is an irreducible $k$-local system of rank $2$ on $\tilde{\PP}^2\setminus \tilde{C},$
having monodromy tuple
$$ \left(   \left(\begin{array}{rr}
1&-2\\
0&1
\end{array}\right),
  \left(\begin{array}{rr}
1&0\\
-2&1
\end{array}\right),
   \left(\begin{array}{rr}
-1&-2\\
2&3
\end{array}\right),
 \left(\begin{array}{rr}
-1&-2\\
2&3
\end{array}\right),
 \left(\begin{array}{rr}
-3&-8\\
2&5
\end{array}\right),
  \left(\begin{array}{rr}
1&-2\\
0&1
\end{array}\right)
\right).$$
\end{prop}

\proof
By the Ogg-Shafarevich Formula~\eqref{degree}, the rank of
$\cR(\cL)$ is $2$ since $n=1,$ $r=4,$ and since $-1\in k^*$ has no nontrivial invariants acting on the
one-dimensional vector space $k.$
Using the theory of wiring diagrams
(encoding the intersection pattern of the underlying real line arrangement, cf.~\cite{Arvola},
\cite{Bartolo}) one derives the braid monodromy of $U$ in a standard way as follows:
$$ \omega_1=\beta_1^2,\, \omega_2={}(\beta_2^2)^{\beta_1},\, \omega_3={}(\beta_1^2)^{\beta_2\beta_1},\,
\omega_4={}(\beta_3^2)^{\beta_1\beta_2\beta_1},\, \omega_5={}(\beta_2^2)^{\beta_3\beta_1\beta_2\beta_1}$$
with $\omega_6$ being the inverse of the product $\omega_1\cdots \omega_5.$

By a computation using
Thm.~\ref{etalem},
 the monodromy tuple of $\cR(\cL)$ is as stated.
 Note that this computation is considerably simplified by the
 obvious observation
 that $\g^{\omega_i}=\g$ for $i=1,\ldots,6.$  Hence  Prop.~\ref{rembraids} in combination with Thm.~\ref{etalem} implies that for the above expression of $\omega_i$ as a word in the standard generators
 $$\omega_i= \beta_{i_1}^{\epsilon_1}\cdots \beta_{i_k}^{\epsilon_k}\quad (\epsilon_j\in \{\pm 1\}),$$ the
matrix inducing $\rho_{\cR(\cL)}(\tilde{\gamma}_i)$ is the corresponding matrix product
 $$M(\g,\omega_i)=M(\g,\beta_{i_1})^{\epsilon_1} \cdots M(\g,\beta_{i_k})^{\epsilon_k}.$$
 The latter product can easily be evaluated using
 $$
M(\g,\beta_1)=\left(\begin{array}{lrlll}
0&\quad -1&\\
1&\quad \,2&&\\
&&&1&\\
&&&&1
\end{array}\right),\,\,
M(\g,\beta_2)=\left(\begin{array}{llrll}
1&&&&\\
&0&\quad -1&\\
&1&\quad \,2&&\\
&&&&1\\
\end{array}\right),\,\,
M(\g,\beta_3)=\left(\begin{array}{lllrl}
1&&&&\\
&1&&&\\
&&0&\quad -1\\
&&1&\quad \,2\\
\end{array}\right)\,
$$
and passing to the subquotient $H_{\g}/E_{\g}.$
 \Endproof

\begin{rem}\label{remstableP}
\begin{enumerate}
\item Again, by standard arguments, the braid monodromy of $ \tilde{C}$
is generated by the following braids:
\begin{equation}\label{eqbraidmon5}  (\beta_4\beta_5\beta_4)^2  ,\quad \beta_3^2,\quad  \beta_3^{-1}(\beta_1\beta_2\beta_1)^2\beta_3
.\end{equation}
It can be checked that the local system $\cR(\cL)$ satisfies indeed the relations coming from
Eq.~\eqref{eqbraidmon5}.
\item
Using the braid monodromy in Eq.~\eqref{eqbraidmon5} the dual Radon transform
applied to   $\cR(\cL)$ has rank $2$ and decomposes into a trivial local system and the restriction of $\cL$
to a slightly smaller open subset. Hence, although the intermediate extension $K=j_{!*}(\cL)[2]$
satisfies the property $P,$ it follows from Prop.~\ref{propP}(2) that
$\cR(K)=\tilde{j}_{!*}\cR(\cL)[2]$ does not  have the property $P.$
\end{enumerate}
\end{rem}

Another example comes from the theory of Zariski pairs (see~\cite{Z2}): Recall the classical Zariski pair,
consisting of a six-cuspidal plane sextic $C$ having its cusps on a conic and a six-cuspidal plane sextic $C'$ whose cusps do not lie on a conic.

Zariski showed that the fundamental group of the complement of $C$ is isomorphic to the free group on two generators of order $2,$ resp. $3,$ whereas the fundamental group of the complement of $C'$ is isomorphic to
a cyclic group of order $6.$ \\

The braid monodromy of the Zariski pair is determined
in \cite{LeeRudolph}, Example 3, (cf.~\cite{Bartolo2}, Prop.~7.3) to be as follows:
for $C$ one has twice the tuple:
$$\left( \beta_1^3, \beta_1^{\beta_2^{-1}\beta_1}, \beta_3^3, \beta_3^{\beta_4^{-1}\beta_3},
\beta_5^3, \beta_2^{\beta_3\beta_3\beta_1},\beta_4^{\beta_5\beta_5\beta_3},\beta_3^{\beta_2^{-1}}, \beta_5^{\beta_4^{-1}}\right).$$
For $C'$ one has:
$$ \left(\beta_3^{\beta_2^{-1}\beta_1\beta_2^{-1}}, \beta_4,\beta_5, \beta_2^3, (\beta_2^{\beta_1^{-1}})^3,\beta_2^{\beta_3^{-1}}, \beta_2^{\beta_1^{-1}\beta_3^{-1}\beta_4^{-1}}, \beta_1^3, \beta_1^{\beta_2^{-1}\beta_3\beta_4^{-1}}, \beta_5^{\beta_4^{-1}\beta_4^{-1}}\right. ,$$
$$\left.
\beta_3^{\beta_2^{-1}\beta_4^{-1}\beta_4^{-1}}, \beta_4^3, \beta_1^3,
\beta_3^{\beta_2^{-1}\beta_1\beta_2^{-1}}, \beta_3^{\beta_2^{-1}\beta_1\beta_2^{-1}},
\beta_1^3, \beta_1^{\beta_2^{-1}\beta_2^{-1}}, \beta_2\right) .
$$
Let now $\cL,$ resp. $\cL',$ be the rank-one $\CC$-local system on $\PP^2\setminus C,$ resp.
$\PP^2\setminus C',$ defined by the obvious homomorphism of fundamental groups
onto the group of 6-th roots of unity inside $\CC^\times=\GL_1(\CC),$ so that in both cases,
the monodromy tuple is a $6$-tuple $(\zeta_6,\ldots, \zeta_6)$ whose entries consist of the same primitive sixth root of unity $\zeta_6\in \CC^\times.$   \\

It is astonishing that the Radon transforms $\cR(\cL),$ resp.  $\cR(\cL'),$
descend to local systems on the complement of the
dual curves $\tilde{\PP}^2\setminus {C}^*,$ resp.
$\tilde{\PP}^2\setminus {C'}^*$ (with  ${C}^*$ being the dual curve of $C$ and
${C'}^*$ being the dual curve of $C'$) since, a priori, they are defined on the complement
of the larger curves $\tilde{C},$ resp. $\tilde{C}',$ defined in Section~\ref{Secradtrans}.  \\

An application of the explicit braid monodromy for the curves $C,$ resp. $C',$  in combination with Thm.~\ref{etalem} leads to the following result:

\begin{thm} \label{thmzariski}
The monodromy group of $\cR(\cL)$ is a finite solvable
group $G\leq \GL_4(\CC)$ of order $648$ acting irreducibly on a direct summand of dimension $3.$ The monodromy group $G'\leq \GL_4(\CC)$ of
$\cR(\cL')$ is a direct product of the cyclic scalar group of order $3$ with the finite symplectic group $\Sp_4(3)\leq \GL_4(\CC).$ \end{thm}

\begin{rem} \begin{enumerate}
\item It is interesting, that the local system $\cL$ does not have the property~$P$ since
its Radon transform splits into a simple subfactor of rank $3$ and a constant subfactor
of rank $1.$

\item In the Appendix to this article, we have listed  the commented
source code which computes the
monodromy tuple of the Radon transform, written in Magma-code~\cite{Magma}.
It can easily be pasted (from the arxiv-latex-source)
to a separate file and then read into Magma.
\end{enumerate}
\end{rem}

\section{Appendix - the Radon program}\label{app:MonMagma}

The following lists the MAGMA code which is used in the proof of Thm.~\ref{thmzariski}.
\begin{small}
\pagestyle{empty}
\begin{lstlisting}

H_g1:=function(g,K)

//Description: Computes the vector space
//$H_{g,1}:=\{(v_1,\ldots, v_r)\mid v_i\in\Im(g_i-I)\}$.

    local r,n,W,I,Kn,i;
    r:=#g;
    n:=Nrows(g[1]);
    W:=VectorSpace(K,0);
    I:=IdentityMatrix(K,n);
    Kn:=VectorSpace(K,n);
    for i in [1..r] do
        W:=DirectSum(W,Image(g[i]-I));
    end for;
    return W;

end function;


H_g2:=function(g,K)

//Description: Computes the vector space
//    $H_{g,2}:=\{(v_1,\ldots, v_r)\mid
//    v_1\cdot g_2\cdot\ldots\cdot g_r+v_2\cdot g_3\cdots g_r+\ldots+v_r=0\}$.

    local r,n,gg,i,A;
    r:=#g;
    n:=Nrows(g[1]);
    gg:=[];
    gg[1]:=g[r];
    A:=ZeroMatrix(K,n*r,n*r);
    for i in [2..r] do
        gg[i]:=g[r-i+1]*gg[i-1];
    end for;
    for i in [1..(r-1)] do
        A:=InsertBlock(A,gg[r-i],(i-1)*n+1,1);
    end for;
    A:=InsertBlock(A,IdentityMatrix(K,n),(r-1)*n+1,1);
    return Kernel(A);

end function;


H_g:=function(g,K)

//Description: Computes the vector space $H_g = H_{g,1}\cap H_{g,2}$

    return (H_g1(g,K) meet H_g2(g,K));

end function;


E_g:=function(g,K)

//Description: Computes the vector space
// $E_g:=\{(v\cdot(g_1-I),\ldots,v\cdot (g_r-I))\mid v\in V\}$.

    local h,r,n,i,I,Im,Mat;
    r:=#g;
    n:=Nrows(g[1]);
    I:=IdentityMatrix(K,n);
    Mat:=ZeroMatrix(K,n,n*r);
    h:=g;
    for i in [1..r] do
        h[i]:=g[i]-I;
    end for;
    Mat:=Hom(VectorSpace(K,n),VectorSpace(K,n*r))!HorizontalJoin(h);
    Im:=Image(Mat);
    return Im;

end function;


Beta_i:=function(g,i)

//Description: Computes the formula for
//the effect of the standard generators $\beta_i$
//  of the Artin braid group $\mathcal{A}_r$ on $\mathcal{E}_r$

//Input:       List $g=[g_1,...,g_r]$, where $g_k \in \GL_n(K)$
//s.t. $\prod_k g_k =1$,
// Integer i, s.t. $|i|\le\#g-1,\,i \neq 0$, representing $\beta_i$

    local r,k,gg;
    r:=#g;
    gg:=g;
    if i ge 1 then
        gg[i]:=g[i+1];
        gg[i+1]:=(g[i+1])^-1*g[i]*g[i+1];
    end if;
    if i le -1 then
        gg[-i]:=g[-i]*g[-i+1]*(g[-i])^-1;
        gg[-i+1]:=g[-i];
    end if;
    return gg;

end function;


phi:=function(g,i,K)

//Description: Computes the map
//$\Phi(g,\beta_i): H_g \rightarrow H_{g^{\beta_i}}$,
//             where $\beta_i$ is a standard generator of $\A_r$
//             (see Section 3).

    local Mat,n,r,gg;
    r:=#g;
    n:=Nrows(g[1]);
    if i ge 1 then
        Mat:=IdentityMatrix(K,n*r);
        Mat:=InsertBlock(Mat,ZeroMatrix(K,n,n),n*(i-1)+1,n*(i-1)+1);
        Mat:=InsertBlock(Mat,IdentityMatrix(K,n),n*i+1,n*(i-1)+1);
            if i le r-1 then
                Mat:=InsertBlock(Mat,g[i+1],n*(i-1)+1,n*i+1);
                Mat:=InsertBlock(Mat,IdentityMatrix(K,n)-(g[i+1])^-1*g[i]*g[i+1],
                               n*i+1,n*i+1);
           end if;
    end if;
    if i le (-1) then
        gg:=Beta_i(g,i);
        Mat:=IdentityMatrix(K,n*r);
        Mat:=InsertBlock(Mat,ZeroMatrix(K,n,n),n*(-i-1)+1,n*(-i-1)+1);
        Mat:=InsertBlock(Mat,IdentityMatrix(K,n),n*-i+1,n*(-i-1)+1);
        if -i le r-1 then
            Mat:=InsertBlock(Mat,gg[-i+1],n*(-i-1)+1,n*(-i)+1);
            Mat:=InsertBlock(Mat,IdentityMatrix(K,n)-(gg[-i+1])^-1*gg[-i]*
                           gg[-i+1],n*(-i)+1,n*(-i)+1);
        end if;
        Mat:=Mat^-1;
    end if;
    return Mat;

end function;


ListenProdukt:=function(L)

//Description: Help function to compute the product over a list
//Input:       List L with entries over a ring
//Output:      Product of the entries of L

    local r,i,LL;
    r:=#L;
    LL:=L[1];
    for i in [2..r] do
        LL:=LL*L[i];
    end for;
    return LL;

end function;


Phi:=function(g,I,K)

//Description: Computes the function $\Phi(g,\beta)
//Input:       List $g=[g_1,...,g_r]$, where $g_k \in \GL_n(K)$
// s.t. $\prod_k g_k =1$,
//  List $I$ of Integers representing the braid $\beta$
//  (e.g. [-1,2,1] represents
// the braid ${\beta_1}^{-1}\cdot\beta_2\cdot\beta_1$
//             Field $K$
//Output:      matrix defining the map $\Phi(g,\beta)$

    local k,L,gg,j;
    k:=#I;
    L:=[];
    gg:=g;
    L[1]:=phi(g,I[1],K);
    for j in [2..k] do
        gg:=Beta_i(gg,I[j-1]);
        L[j]:=phi(gg,I[j],K);
    end for;
    return ListenProdukt(L);

end function;


Trafodat:=function(g,K)

    local B1,B2,B3,M,n,r;
    r:=#g;
    n:=Nrows(g[1]);
    B1:=Basis(E_g(g,K));
    B2:=ExtendBasis(B1,H_g(g,K));
    B3:=ExtendBasis(B2,VectorSpace(K,n*r));
    M:=Matrix(K,n*r,n*r,B3);
    return [*Dimension(E_g(g,K)),Dimension(H_g(g,K)),M*];

end function;


Phibar:=function(g,I,K)

//Description: Computes the isomorphism
//             $\bar{\Phi}(g,\beta):
//H_g/E_g \rightarrow H_{g^\beta}/E_{g^\beta}$,
//             induced by the map $\Phi(g,\beta)$
//Input:       List $g=[g_1,...,g_r]$, where
// $g_k \in \GL_n(K)$ s.t. $\prod_k g_k =1$,
//  List $I$ of Integers representing the braid $\beta$
//             Field $K$
//Output:      matrix, representing the
// isomorphism $\bar{\Phi}(g,\beta)$

    local n,r,P,T,M;
    r:=#g;
    n:=Nrows(g[1]);
    P:=Phi(g,I,K);
    T:=Trafodat(g,K);
    M:=T[3]*P*(T[3])^-1;
    return ExtractBlock(M,T[1]+1,T[1]+1,T[2]-T[1],T[2]-T[1]);

end function;


Radon:=function(g,W,K)

//Computes the monodromy tuple of the Radon transform
// Input: $g=[g_1,...,g_r]$, where
//$g_k \in \GL_n(K)$ s.t. $\prod_k g_k =1$
//$W=[W[1],...,W[s]] \in \cB_r$
//braid monodromy of the singular locus $C$
// of the local system belonging to $g$.

local n,r,s,TT,i;
    r:=#g;
    s:=#W;
    n:=Nrows(g[1]);
    TT:=[];

    for i in [1..s] do
         TT[i]:=Phibar(g,W[i],K);
    end for;

return TT;

end function;




\end{lstlisting}
\end{small}
\end{document}